\documentclass{amsart}
\numberwithin{equation}{section}
\newcommand{\cH}{\mathcal{H}}
\newcommand{\cR}{\mathcal{R}}
\newcommand{\cW}{\mathcal{W}}
\newcommand{\cC}{\mathcal{C}}
\newcommand{\cB}{\mathcal{B}}
\newcommand{\cV}{\mathcal{V}}
\newcommand{\bL}{\mathbf{L}}
\newcommand{\CC}{\mathbb{C}}
\newcommand{\DD}{\mathbb{D}}
\newcommand{\MM}{\mathbb{M}}
\newtheorem{Pa}{Paper}[section]
\newtheorem{Tm}[Pa]{{\bf Theorem}}

\newtheorem{Rk}[Pa]{{\bf Remark}}
\newtheorem{Pn}[Pa]{{\bf Proposition}}
\newtheorem{Pb}[Pa]{{\bf Problem}}

\DeclareMathOperator{\ran}{ran}
\DeclareMathOperator{\diag}{diag}

\begin{document}
\begin{abstract}
We study Carath\'eodory-Herglotz functions whose values are 
continuous operators
from a locally convex topological space which admits the
factorization property into its conjugate dual space. We show how
this case can be reduced to the case of functions whose values
are bounded operators from a Hilbert space into itself.
\end{abstract}
\title[Carath\'eodory--Fej\'er interpolation]
{Carath\'eodory--Fej\'er interpolation and related topics in
locally convex spaces}

\author{Daniel Alpay, Olga Timoshenko and Dan Volok}



\date{}

\maketitle

\section{Introduction}

In the present paper we pursue our study of Carath\'eodory-Herglotz functions in
the case of operator valued functions. In \cite{atv1} we studied
Carath\'eodory-Herglotz\footnote{We used the terminology {\sl Carath\'eodory
function} in \cite{atv1}.} functions whose values are bounded
operators from a Banach space $\cB$ into its conjugate
dual space $\cB^*$ (that is, into the space of its
anti-linear continuous functionals). The present work is set in
the framework of locally convex vector spaces which have the
factorization property. The methods and focus in \cite{atv1} are
different than those in the present paper. In \cite{atv1} we used
in an essential way the theory of linear relations in Hilbert
spaces, and we obtained realization theorems for Carath\'eodory-Herglotz
functions; see Theorem \ref{michelle} below. Here the
factorization property for positive operators in certain locally
convex topological vector spaces plays a central role, and the
main point is the reduction to the Hilbert space case; see Theorem
\ref{montpellier} below.\\

Let  $\cV$ be a locally
convex topological vector space, let $\cV^*$ denote
the space of its antilinear continuous functionals and let
\[\langle\cdot,\cdot\rangle_{\cV}\]
denote the duality between $\cV$ and
$\cV^*.$ We shall denote by
$\bL(\cV,\cV^*)$
the space of continuous operators from $\cV$ into $\cV^*$ and, more generally,
by $\bL(\cV,\cW)$ the space of continuous operators from $\cV$ into another locally
convex topological vector space $\cW.$ For a Hilbert space $\cH$, $\bL(\cH)$ stands for the space of bounded operators from $\cH$ into itself. \\

An operator $A\in\bL(\cV,\cV^*)$ is said to be positive  (notation: $A\geq 0$) if
\begin{equation}
\langle Av,v\rangle_{\mathcal{V}}\ge 0,\quad
\forall v\in\cV.
\end{equation}
The space $\cV$ is said to possess the  factorization property if
every  $A\geq 0$ in $\bL(\cV,\cV^*)$ admits a factorization of the form
\[A=T^*T,\]
where $T$ is a continuous operator from $\cV$ into some Hilbert space (depending on $A$).
It was proved by
J. Gorniak and A. Weron (see \cite[Theorem 2.13 p. 240]{MR647140}) that $\cV$ has the  factorization property
if and only if for every $A\geq 0$ in $\bL(\cV,\cV^*)$ the map
\[
v\mapsto \langle Av,v\rangle_{\mathcal{V}}
\]
is continuous. We remark that the results of \cite{atv1} are in fact valid in any
locally convex vector space which has the factorization
property, and this is the setting in which we
consider the present work.\\

An $\bL(\cV,\cV^*)$-valued function
$\Phi$ defined in the open unit disk $\DD$ (but {\it a
priori} without any analyticity property) is called a
Carath\'eodory-Herglotz function if the kernel
\begin{equation}
\label{190607}
K_\Phi(z,w)=\dfrac{\Phi(z)+\Phi(w)^*\big|_{\mathcal
V}}{1-z\overline{w}}
\end{equation}
is positive in $\DD$. We denote by $\cC(\cV, \cV^*)$ the family of $\bL(\cV,\cV^*)$-valued Carath\'eodory-Herglotz functions. In the case when $\cV$ is a Hilbert space $\cH$ we use the
notation $\cC(\cH)$.\\

\begin{Tm}
\label{michelle} Let $\cV$ be a locally convex
topological vector space with the factorization property,
and let $\Phi\in \cC(\cV, \cV^*)$. Then
there exist a Hilbert space $\cH$ and operators $C\in\bL(\cV,\cH),$
$V\in\bL(\cH),$
$D\in\bL(\cV,\cV^*),$ such that:
\begin{enumerate}
\item
$V$ is an isometry: $$V^*V=I_\cH$$ and $D$ is purely imaginary:
\[
D+D^*\big|_{\cV}=0;\]
\item $\Phi(z)$ admits the representation
\begin{equation}
\label{sarah}
\Phi(z)=D+C^*(I+zV^*)(I-zV^*)^{-1}C,\quad z\in{\mathbb D}.
\end{equation}
In particular, $\Phi(z)$ admits (in the weak sense) the power series expansion
\begin{equation}
\label{panthere}
\Phi(z)=M_0+2\sum_{n=1}^\infty z^nM_n,
\end{equation}
where the coefficients $M_n\in\bL(\cV,\cV^*)$ are given by
\begin{equation}
\label{oasis}
M_n=\begin{cases} D+C^*C,&\quad \text{if}\quad n=0,\\
C^*V^{*n}C,&\quad \text{if}\quad n\geq 1,\end{cases}
\end{equation}
and, therefore,
 $\Phi(z)$ is weakly analytic\footnote{See \cite[p. 80-81]{schaefer}.} in the open unit disk $\DD$.
 \end{enumerate}
\end{Tm}
\begin{proof}
The representation \eqref{sarah} was originally obtained in \cite{atv1} in the setting of
Banach spaces;
the proof in the present setting is absolutely the same.
 \end{proof}
\begin{Rk}In the case when $\cV$ is a Banach space,
 $\bL(\cV,\cV^*)$ is a Banach space, as well, and
weak analyticity is equivalent to strong analyticity; see
\cite[Theorem VI.4, p. 189]{MR58:12429c}. Then $\Phi$
is strongly analytic, which is also evident
 from the fact that the coefficients $M_n$ in \eqref{panthere} are uniformly
bounded.
\end{Rk}

\begin{Rk}Formulas \eqref{oasis} are the operator-theoretic version of
formulas appearing in the theory of stationary stochastic
processes; see for instance \cite[p. 45]{MR0348830}.\end{Rk}

 From \eqref{oasis} we see that the operator matrices
\begin{equation}\label{16}
\begin{split}
&\begin{pmatrix}
\frac{M_0+M_0^*|_\cV}{2}&M_1^*|_\cV&M_2^*|_\cV
&\cdots &M_N^*|_\cV\\
M_1 &\frac{M_0+M_0^*|_{\mathcal
V}}{2}&M_1^*|_\cV&\cdots&M_{N-1}^*|_\cV\\
M_2&M_1&\frac{M_0+M_0^*|_\cV}{2}&\cdots&\\
\vdots&\cdots& & &\\
M_N&\cdot& & &\frac{M_0+M_0^*|_\cV}{2}&
\end{pmatrix}\\[1ex]
&\qquad =\diag(C^*,\ldots, C^*)
\begin{pmatrix}I_\cH\\ V^*\\ \vdots \\
  V^{*N}
\end{pmatrix}
\begin{pmatrix}I_\cH& V& \cdots &
  V^{N}
\end{pmatrix}
\diag(C,\ldots, C),
\end{split}\end{equation}
where $N=0,1,\ldots$, are positive elements of $\bL(\cV^{N+1},(\cV^*)^{N+1})$.
Post- and premultiplying \eqref{16} with the permutation matrix
\[\begin{pmatrix} 0&\cdots&0&I_\cV\\
0&&I_\cV&0\\
\hdotsfor[2]{4}\\
I_\cV&0&\cdots&0\end{pmatrix}\]
and its adjoint, respectively, we conclude that
 the operator matrices
\begin{equation}
\label{natasha}
\MM_N=\begin{pmatrix}
\frac{M_0+M_0^*|_\cV}{2}&M_1&M_2&\cdots &M_N\\
M_1^*|_\cV&\frac{M_0+M_0^*|_{\mathcal
V}}{2}&M_1&\cdots&M_{N-1}\\
M_2^*|_\cV&M_1^*|_{\mathcal
V}&\frac{M_0+M_0^*|_\cV}{2}&\cdots&\\
\vdots&\cdots& & &\\
M_N^*|_\cV&\cdot& & &\frac{M_0+M_0^*|_\cV}{2}
\end{pmatrix},
\end{equation}
are positive, as well.\\

As we shall see in the sequel (Theorem \ref{montpellier} below), the
converse is also true: if all the matrices $\MM_N$ are
positive, the corresponding power series \eqref{panthere} defines an
element in $\cC(\cV,\cV^*)$.\\

We pose
the following Carath\'eodory-Fej\'er  interpolation problem:

\begin{Pb}
\label{CF}
Given $N\in{\mathbb N}$ and operators $M_0,\ldots
M_{N}\in\bL(\cV,\cV^*)$, such that the Toeplitz operator matrix $\MM_N$ is positive,
find a
 function
$\Phi\in\cC(\cV,\cV^*)$ such
that
\[
\Phi(z)-(M_0+2\sum_{n=1}^{N}z^nM_n)=O(z^{N+1}).
\]
\end{Pb}

We shall prove (see Theorem \ref{Pont Marie} below)  that Problem \ref{CF} is always solvable.
Furthermore, it turns out that
 $\Phi\in\cC({\mathcal
V},\cV^*)$  if and only if
\begin{equation}
\label{ertyu}
\Phi(z)=D+T^*\varphi(z)T,
\end{equation}
where
$D\in\bL(\cV,\cV^*)$ is purely imaginary,
 $T\in\bL(\cV,\cH)$ for some Hilbert space
$\cH$ and $\overline{\ran(T)}=\cH,$  $\varphi\in{\mathcal
C}(\cH)$. Moreover, given $\Phi,$ the operator $D$ is determined uniquely while the operator
$T$ and the space $\cH$ are determined up to unitary mappings.
 Hence we can reduce interpolation
problems for which the value of the function $\Phi\in{\mathcal
C}(\cV,\cV^*)$ is pre-assigned at the origin to
interpolation problems for functions
$\varphi\in\cC(\cH)$.\\

This last point explains also why formula \eqref{ertyu} is much
more precise and useful than the seemingly similar formula
\eqref{sarah}. In \eqref{sarah} the operator $C$ need not have
dense range (in fact $C$ can be chosen to be the adjoint of point
evaluation at the origin), and the reduction
mentioned above does not apply for \eqref{sarah}.\\

The original motivation for \cite{atv1} and for the present work
originates with works on extrapolation stationary stochastic
processes with values in a Banach space, or more generally in a
linear space; see for instance \cite{MR611719}, \cite{MR521020},
\cite{MR521027}. Our aim is to develop the theory of de
Branges-Rovnyak spaces, and associated problems such as inverse
scattering (see \cite{ad1}, \cite{ad2}) in the linear space
case.\\

 The paper consists of four sections besides the introduction.
In the second section we review some facts on locally vector
spaces which have the factorization property. In the third section 
we review the solution of the Carath\'eodory-Fej\'er in the Hilbert
space case. The main results of this paper, including the
characterization of a Carath\'eodory-Herglotz function in terms of its
power series and the solution of the
Carath\'eodory-Fej\'er  interpolation problem appear in Section 4.
\section{Preliminaries: spaces with the factorization property}
Let $\cV$ be a locally convex
topological vector space with the factorization property and let $A\in\bL(\cV,\cV^*)$ be positive.
We shall say that the factorization
\(A=T^*T,\)
where $T\in\bL(\cV,\cH)$ and $\cH$ is a Hilbert space, is {\em minimal} if
the range of $T$ is dense in $\cH$.
\begin{Pn}
\label{factor} Let $\cV$ be a locally convex
topological vector space with the factorization property and let $A\in\bL(\cV,\cV^*)$ be positive. Let
\begin{align}\label{f1}A&=T^*T,\\ \label{f2}A& =(T^{\prime})^*T^\prime\end{align}
be two factorizations via Hilbert spaces $\cH$ and
 $\cH^\prime$ and assume that the factorization \eqref{f1} is minimal.  Then there exists an isometry
$V\in\bL(\cH,\cH^\prime)$ such that $T^\prime =VT$. Moreover, if the factorization \eqref{f2} is also minimal then
$V$ is unitary.
\end{Pn}
\begin{proof}
The proof follows the argument in \cite{MR647140}. The formula
\[
\cR=\left\{(Tv, T^\prime v),\quad v\in\cV\right\}.
\]
define a relation on $\cH\times\cH^\prime$, which is
isometric and has a dense domain; it is therefore the graph of an
isometric operator $V$. If the range of $\cR$ is dense, as well, then $V$ must be unitary.
\end{proof}

\begin{Pn}
\label{Invalides}
Assume that $\cV$ admits the
factorization property. Then so does $\cV^N$.
\end{Pn}

\begin{proof} Let $A$ be a positive continuous operator from $\cV^N$
into $(\cV^*)^N$. We consider the matrix representation
$A=(A_{\ell j})_{\ell, j=1,\ldots N}$, where $A_{\ell j}\in\bL(\cV,\cV^*)$.
We want to prove that the map
\begin{equation}
(v_1,\ldots v_n)\mapsto\left\langle A\begin{pmatrix}v_1\\
\vdots \\v_N\end{pmatrix},\begin{pmatrix}v_1\\ \vdots\\
v_N\end{pmatrix}\right\rangle_{\cV^N}
\label{yesterday}
\end{equation}
is continuous. The positivity of $A$ implies that the operators
$A_{\ell\ell}$ are positive. Since $\cV$ has the
factorization property, the maps
\[
v\mapsto \langle A_{\ell \ell} v,v\rangle_{\cV},\quad \ell=1,\ldots N,
\]
are continuous. Thus it suffices to show that the maps
\begin{equation}
(v,w)\mapsto \langle A_{\ell j} w,v\rangle_{\cV} \label{perlette}
\end{equation}
are continuous for $\ell\not= j$.\\

Without loss of generality, we assume that $\ell<j$. The positivity of $A$ implies that
the operator
\[
\begin{pmatrix}
A_{\ell\ell}&A_{\ell j}\\
A{j\ell}&A_{jj}\end{pmatrix}
\]
is positive from $\cV^2$ in to $(\cV^*)^2$. It
follows that for every choice of $v,w\in\cV$ and
$\alpha,\beta\in\CC$
\[
\begin{pmatrix}\overline{\alpha}&\overline{\beta}\end{pmatrix}
\begin{pmatrix}\langle A_{\ell\ell}v,v\rangle_{\cV} & \langle A_{\ell j}w,v\rangle_{\cV}
\\
\langle A_{j\ell}v,w\rangle_{\cV} &
\langle A_{jj}w,w\rangle_{\cV}\end{pmatrix}\begin{pmatrix}\alpha\\
\beta\end{pmatrix}\ge 0.
\]
This inequality implies first that
\[
\langle A_{j\ell}v,w\rangle_{\cV}=
\overline{ \langle A_{\ell j}w,v\rangle_{\cV}},\]
and so
\begin{equation}
A_{\ell j}^*\big|_\cV=A_{j\ell}.
\end{equation}
Moreover the matrix
\[
\begin{pmatrix}\langle A_{\ell\ell}v,v\rangle_{\cV} & \langle A_{\ell j}w,v\rangle_{\cV}\\ \langle A_{j\ell}v,w\rangle_{\cV} &
\langle AA_{jj}w,w\rangle_{\cV}\end{pmatrix}
\]
is positive, and hence so is its determinant, that is:
\[
|\langle A_{\ell j}w,v\rangle_{\cV}|^2\le
\langle A_{\ell\ell}v,v\rangle_{\cV}
\langle A_{jj}w,w\rangle_{\cV}.\]
We deduce that the map \eqref{perlette} is continuous. The
continuity of the map \eqref{yesterday} follows.
\end{proof}

\section{The Carath\'eodory--Fej\'er extension
problem in the Hilbert space case}
\setcounter{equation}{0}
We recall a result on the Carath\'eodory--Fej\'er extension
problem in the Hilbert space case. We give a proof for
completeness, using a one step extension procedure, as in for
instance \cite{van-den-bos}. Note that the
spaces are not assumed to be separable. \\

We shall use the following formula (see \cite[(0.3) p. 3]{MR90g:47003}):
let \[
G=\begin{pmatrix}A&B\\ C&D\end{pmatrix}
\]
be
such that $A$ is invertible, then
\begin{equation}
\label{useful formula}
G=\begin{pmatrix} I&0\\ CA^{-1}& I\end{pmatrix}\begin{pmatrix}
A&0\\0&D^\times\end{pmatrix}\begin{pmatrix}I&A^{-1}B\\0&I
\end{pmatrix},\end{equation}
where $D^{\times}=D-CA^{-1}B$ is called the Schur complement of $D$.

\begin{Tm}
If $\cV=\cH$ is a Hilbert space then Problem \ref{CF} is solvable.
\label{Wagram, ligne 3}
\end{Tm}
\begin{proof} We replace $\MM_N$ by $\epsilon I_{\cH^{N +1}}+\MM_N$, where $\epsilon>0$. This last operator
is strictly positive (that is, positive and boundedly invertible): 
\[\epsilon I_{\cH^{N +1}}+\MM_N>0.\]
 We set
$$(\epsilon I_{\cH^{N+1}}+
 \MM_N)^{-1}=\begin{pmatrix} \alpha_N & \beta_N^*\\
\beta_N & \delta_N
\end{pmatrix}\quad\text{and}\quad\begin{pmatrix}
M_N &\ldots & M_1\end{pmatrix}=\gamma_N,$$ where, to lighten
the notation, we do not stress for the moment the dependence
 on $\epsilon$. Let $X\in\bL(\cH)$. In view of \eqref{useful formula}, we see that the Toeplitz matrix
\[
\begin{pmatrix}
\epsilon I_{\cH^{N+1}}+\MM_N&\begin{array}{c}X^*\\\gamma_N^*\end{array} \\
X\qquad\gamma_N &\epsilon I_{\cH}+M_0\end{pmatrix}\]
is strictly positive if and only if
\[
\epsilon I_\cH+M_0>\begin{pmatrix} X & \gamma_N \end{pmatrix} (\epsilon I_{\cH^{N+1}}+\MM_N)^{-1} \begin{pmatrix}
X^*\\
\gamma_N^*\end{pmatrix}.
\]
Using \ref{useful formula}, the last inequality can be rewritten as 
\begin{equation}
\label{ball}\epsilon I_{\mathcal H}+M_0-\gamma_N\delta_N^\times\gamma_N^*>(X+\gamma_N\beta_N\alpha_N^{-1})\alpha_N(X+\gamma_N\beta_N\alpha_N^{-1})^*.
\end{equation}
But the operator on the left hand side of \eqref{ball} is the Schur complement of $\epsilon
I_\cH+M_0$ in $\epsilon I_{\cH^{N +1}}+\MM_N,$ which is strictly positive. Hence
 the operator ball defined by \eqref{ball} is non-empty.
 Setting $M_{N+1}(\epsilon)=X,$ where
  $X$
satisfies \eqref{ball}, we obtain 
$\MM_{N+1}>0.$\\

Reiterating this procedure we can choose a sequence of operators
\[M_{N+1}(\epsilon),M_{N+2}(\epsilon),\ldots\]
such that all the
corresponding Toeplitz matrices
\[\MM_{N+1}(\epsilon),\MM_{N+2}(\epsilon),\ldots\]
are strictly positive. The corresponding function
\[
\Phi_\epsilon(z)=M_0+2\sum_{k=1}^N M_k z^k+2\sum_{k=N+1}^\infty
M_k(\epsilon) z^k
\]
is in the Carath\'eodory-Herglotz class since
\begin{multline}
\label{Guy}
\frac{\Phi_\epsilon(z)+\Phi_\epsilon(w)^*}{1-z\overline{w}}\\=2\lim_{n\rightarrow\infty}
\begin{pmatrix} z^nI_\cH&z^{n-1}I_\cH&\cdots &I_\cH
\end{pmatrix}\MM_{n}(\epsilon)
\begin{pmatrix}{w}^n I_\cH&{w}^{n-1}I_\cH&\cdots& I_\cH
\end{pmatrix}^*.
\end{multline}

The positivity of the matrices $\MM_{n}(X)$ implies, much
in the same way as in the proof of Proposition \ref{Invalides},
that the operators $M_n(\epsilon)$ are uniformly bounded in norm
for $\epsilon\le 1$. We now want to let $\epsilon\rightarrow 0$.
We use the fact (see for instance \cite[Chapter 5]{beauzamy})
that $\bL(\cH)$ is the dual space of the trace
class operators in $\bL(\cH)$, and hence,
by Banach-Alaoglu theorem (see \cite[Lemma 2.3, p.
102]{beauzamy})
 the closed unit ball of the space
$\bL(\cH)$ is weakly compact.  Let now $\epsilon_k$,
$k=0,1,\ldots$ be a sequence of numbers decreasing to $0$. Since
rank one operators are trace class, we can find a subsequence
$\epsilon_{1,k}$ and an operator $M_{N+1}\in\bL(\cH)$ such that
\[
\lim_{k\rightarrow \infty}\langle
M_{N+1}(\epsilon_{1,k})x,y\rangle_\cH= \langle
M_{N+1}x,y\rangle_\cH,\quad\forall x,y\in\cH.
\]
Similarly there exists a subsequence  $\epsilon_{2,k}$ of the
sequence $\epsilon_{1,k}$ and an operator $M_{N+2}\in\bL(\cH)$ such
that
\[
\lim_{k\rightarrow \infty}\langle
M_{N+2}(\epsilon_{2,k})x,y\rangle_\cH= \langle
M_{N+2}x,y\rangle_\cH,\quad\forall x,y\in\cH.
\]
Using the diagonal argument, we find   a sequence of bounded
operators $M_n,n>N$ and a sequence of numbers $\eta_k=\epsilon_{k,k}$ decreasing
to $0$ such that
\[
\lim_{k\rightarrow\infty}\langle \Phi_{\eta_k}(z)v,w\rangle=
\langle \Phi(z)v,w\rangle,\]
where \( \Phi(z)=M_0+2\sum_{n=1}^{\infty}z^nM_n.
\) The function $\Phi$ is a Carath\'eodory-Herglotz function. Indeed,
for every $k$ the kernel $K_{\Phi_{\eta_k}}(z,w)$ is positive.
Thus for every choice of $n\in{\mathbb N}$, $z_1,\ldots,
z_n\in{\mathbb D}$ and $h_1,\ldots, h_n\in\cH$ the
$n\times n$ matrix with $\ell j$ entry equal to
\[
\langle K_{\Phi_{\eta_k}}(z_\ell,z_j)h_j,h_\ell\rangle_{\mathcal
H}\]
is positive. Letting $k\rightarrow\infty$ we get that the kernel
$K_{\Phi}(z,w)$ is positive in ${\mathbb D}$.
\end{proof}\

\section{Characterization of a Carath\'eodory-Herglotz function in terms of its
power series}
As already noticed, the space $\cV^N$ has the
factorization property when $\mathcal V$ has. As a consequence we
have:

\begin{Tm}
\label{Pont Marie}
Let $N\in{\mathbb N}$ and let $M_0,\ldots,
M_N$ be continuous operators from $\mathcal{V}$ into $\mathcal{V}^*.$ Assume that the
Toeplitz block operator matrix $\MM_N$ given by
\eqref{natasha} is a positive operator from $\cV^{N+1}$ into
$(\cV^*)^{N+1}$.
Let
\[
\frac{M_0+M_0^*|_\cV}{2}=T_0^*T_0,
\]
where $T_0$ is a continuous operator from $\mathcal{V}$
 into a Hilbert space $\cH_0$ with dense range.
 Then there exist a unitary operator $U$ from $\cH_0$ into
itself and a bounded operator $C$ from $\cH_0$ into
itself such that
\begin{equation}
\frac{M_0+M_0^*|_\cV}{2}=T_0^*C^*CT_0\quad{and}\quad
M_j=T_0^*C^*U^jCT_0,\quad,j=1,\ldots N.
\end{equation}
In particular, the function
\[
\Phi(z)=\frac{M_0-M_0^*\big|_\cV}{2}+T_0^*
C^*(I+zU)(I-zU)^{-1}CT_0
\]
is a solution of the corresponding Carath\'eodory--Fej\'er extension
problem.
\end{Tm}
{\bf Proof:} Since $\cV^{N+1}$ has the factorization
property, there exists an Hilbert space $\cH$ and a
continuous operator $T$ from $\cV^{N+1}$ into ${\mathcal
H}$ such that $\MM_N=F^*F.$ Write
\[
F=\begin{pmatrix}F_0&F_1&\cdots&F_N\end{pmatrix}.
\]
We have in particular
\[
\frac{M_0+M_0^*|_\cV}{2}=F_0^*F_0=F_1^*F_1=\cdots=F_N^*F_N.
\]
By Proposition \ref{factor},  there exist isometries
$V_0,\ldots, V_N$ from $\cH_0$ into $\cH$ such
that
\[
F_j=V_jT_0,\quad j=0,\ldots, N.\]
Let $j>i$. The equality $\MM_N=F^*F$ leads to
\begin{equation}
\label{Boucicaut, ligne 8}
M_j=T^*_0V_0^*V_jT_0=T^*_0V_1^*V_{j-1}T_0=\cdots
=T_0^*V_i^*V_{j-i}T_0.
\end{equation}
In view of \eqref{Boucicaut, ligne 8} we have:
\begin{equation}
\label{Lourmel, ligne 8}
\begin{split}
\MM_N&=\begin{pmatrix}T_0^* &0    &0&\cdots &0\\
                             0    &T_0^*&0&\cdots&0\\
                             \cdot&\cdot&\cdot&\cdot\\
                              \cdot&\cdot&\cdot&\cdot\\
                             0&\cdot& &
                             &T_0^*\end{pmatrix}\times\\
&\hspace{5mm}\times
\begin{pmatrix}I&V_0^*V_1&\cdots&V_0^*V_N\\
               V_1^*V_0&I&  V_0^*V_1&    \cdots & \\
               V_2^*V_0&V_1^*V_0&I&\cdot\\
                \cdot&\cdot&\cdot&\cdot\\
               V_N^*V_0& & & I
               \end{pmatrix}
                             \begin{pmatrix}T_0 &0    &0&\cdots &0\\
                             0    &T_0&0&\cdots&0\\
                             \cdot&\cdot&\cdot&\cdot\\
                              \cdot&\cdot&\cdot&\cdot\\
                             0&\cdot& & &T_0\end{pmatrix}.
\end{split}
                             \end{equation}
Furthermore the block Toeplitz operator
\[
\begin{pmatrix}I&V_0^*V_1&\cdots&V_0^*V_N\\
               V_1^*V_0&I&  V_0^*V_1&    \cdots & \\
               \cdot&\cdot&\cdot&\cdot\\
                \cdot&\cdot&\cdot&\cdot\\
               V_N^*V_0& & & I
               \end{pmatrix}\ge 0
               \]
since $T_0$ has dense range. The problem is thus reduced to the
Hilbert space case, and Theorem \ref{Wagram, ligne 3} is
applicable. The proof is now easily completed.
\mbox{}\qed\mbox{}

We now turn to the main result of this paper:
\begin{Tm}\label{montpellier}
Let $M_0,M_1,\ldots$ be an infinite sequence of elements of
$\bL(\cV,\cV^*)$ such that all the
matrices $\MM_N$ defined by \eqref{natasha},
$N=0,1,\ldots$, are positive. Then for every $z$ in the open
unit disk, the series
\[
\Phi(z)=M_0+2\sum_{n=1}^\infty z^nM_n
\]
converges weakly, and defines a Carath\'eodory-Herglotz function. Furthermore,
$\Phi$ can be written as \eqref{ertyu}:
\begin{equation*}
\Phi(z)=\dfrac{\Phi(0)-\Phi(0)^*\big|_{\mathcal
Vv}}{2i}+T_0^*\varphi(z)T_0.
\end{equation*}
Conversely, every function of the form \eqref{ertyu} belongs to 
$\cC({\mathcal V},\cV^*)$.
\end{Tm}
{\bf Proof:} From the proof of Theorem \ref{Pont Marie} follows
that that there is a Hilbert space $\cH_0$ and an
operator $T_0\in\bL({\mathcal B}, \cH_0)$ with
dense range and such that for every $N\ge 0$,
\[
{\mathbf M}_N={\rm diag}~(T_0^*,\ldots T_0^*){\mathbf T}_N{\rm
diag}~(T_0,\ldots ,T_0),
\]
where ${\mathbf T}_N$ is a uniquely defined positive block
Toeplitz operator from $\cH_0^{N+1}$ into ${\mathcal
H}_0^{N+1}$. Since the ${\mathbf T}_N$ are uniquely defined, they
are of the form
\[
{\mathbf T}_N=(t_{i-j})_{i,j=0,\ldots N},\]
where the $t_j$ are bounded operators from $\cH_0$ into
itself. As in the proof of Theorem \ref{Wagram, ligne 3} (see
\eqref{Guy}),
the function
\[
\varphi(z)=t_0+2\sum_{n=1}^\infty t_nz^n
\]
belongs to ${\mathcal C}({\mathcal H}_0)$, and \eqref{ertyu} holds
with $\varphi$.\\

The converse follows from
\[
K_\Phi(z,w)=T_0^*K_\varphi(z,w)T_0.
\]
\mbox{}\qed\mbox{}\\

\bibliographystyle{plain}
\def\cprime{$'$} \def\cprime{$'$} \def\cprime{$'$} \def\cprime{$'$}
  \def\cprime{$'$}

\end{document}